\documentclass[11pt]{amsart}
\usepackage[margin=2.5cm]{geometry}
\usepackage{amssymb}
\usepackage{hyperref} 
\usepackage{bm}
\usepackage[centertags]{amsmath}
\usepackage{amsfonts}
\usepackage{amsthm}
\usepackage{mathrsfs}
\usepackage{mathtools,xfrac}
\usepackage[usenames,dvipsnames,svgnames,table]{xcolor}
\usepackage[mathscr]{euscript}
\usepackage[shortlabels]{enumitem}

\DeclareMathOperator\w{w}
\DeclareMathOperator\range{rng}

\DeclareMathOperator\supp{supp}

\DeclareMathOperator\suc{succ}

\DeclareMathOperator{\spa}{span}

\newcommand{\mt}{\mathcal{T}}

\numberwithin{equation}{section}

\newtheorem{theorem}{Theorem}[section]
\newtheorem*{theorem*}{Theorem}
\newtheorem{lemma}[theorem]{Lemma}
\newtheorem{proposition}[theorem]{Proposition}
\newtheorem{corollary}[theorem]{Corollary}

\newtheorem{fact}[theorem]{Fact}

\newtheorem{definition}[theorem]{Definition}

\newtheorem{notation}[theorem]{Notation}
\newtheorem{assumption}[theorem]{Assumption}

\newtheorem{remark}[theorem]{Remark}

\newcommand{\N}{\mathbb{N}}
\newcommand{\Z}{\mathbb{Z}}

\newcommand{\R}{\mathbb{R}}
\newcommand{\Q}{\mathbb{Q}}

\usepackage[utf8]{inputenc}

\author{Anna Pelczar-Barwacz}
\thanks{The research of the author was supported by the grant of the National Science Centre (NCN), Poland, no. UMO-2020/39/B/ST1/01042}

\subjclass[2020]{46B04, 46B15, 46B20, 46B45}

\email{anna.pelczar@uj.edu.pl}
\address{Jagiellonian University, Faculty of Mathematics and Computer Science, Institute of Mathematics,  {\L}ojasiewicza 6, 30-348 Krak\'ow,  Poland}

\usepackage[T1]{fontenc}
\numberwithin{subsection}{section}
\numberwithin{equation}{section}

\begin{document}
\title{An example of a reflexive Banach space with an uncomplemented isometric subspace} 

\begin{abstract}
We present an example of a reflexive Banach space $X$ with a  basis $(e_i)_{i=1}^\infty$ 1-equivalent to its subsequence $(e_{2i})_{i=1}^\infty$. The closed subspace of $X$ spanned by the latter 
is not complemented in $X$, whereas the subsequence $(e_{2i-1})_{i=1}^\infty$ is unconditional. 
\end{abstract}
\maketitle

\section{Introduction}

The study of isometries on Banach spaces usually concentrates on the group of surjective isometries. The natural question in the research on partial  isometries concerns the complementability of the range of an isometry in the domain space. The strongest version of the latter property, i.e. 1-complementability of the range of an isometry, allows in the case of reflexive Banach spaces for an extended version of the Wold decomposition with respect to the isometry in question, \cite{cfs}, well known in the case of Hilbert spaces. The class of Banach spaces, on which all isometries have  1-complemented ranges includes Hilbert spaces, $L_p$-spaces, $1\leq p<\infty$, on any measure space \cite{bl,pelczynski,tzafriri}. On the other hand, $C[0,1]$ contains an isometric subspace that is not complemented \cite{ditor}; whereas in the reflexive range, each $\ell_p$, $2<p<\infty$ contains an uncomplemented  isomorphic subspace $Y_p$ by \cite{rosenthal}. As mentioned above, such a subspace cannot be isometric to $\ell_p$, but can be used to built a reflexive space with uncomplemented isometric subspace: let $W_p:=Y_p\oplus_p \ell_p\subset \ell_p\oplus\ell_p\equiv \ell_p\subset Y_p\oplus\ell_p=: Z_p$, with all inclusions isometric. It follows that $W_p$ is not complemented in $Z_p$, as otherwise $Y_p$ would be complemented in $\ell_p$. We thank Tomasz Kania for the possibility of including here his example above. For other results described above see an excellent survey \cite{randria}.

We present another example of a reflexive Banach space with an uncomplemented isometric subspace within the framework of basic sequences, namely we built a reflexive Banach space with a basis, that is equivalent to some of its subsequences  spanning an uncomplemented subspace. Obviously, the constructed basis cannot be unconditional itself, however, it contains an unconditional subsequence. The precise result reads as follows. 
\begin{theorem}\label{main}
There is a reflexive Banach space $X$ with a basis $(e_i)_{i=1}^\infty$ satisfying the following.
\begin{enumerate}[(1)]
\item\label{main:2} Sequences $(e_i)_{i=1}^\infty$ and $(e_{2i})_{i=1}^\infty$ are 1-equivalent, that is
the map carrying each  $e_i$ to $e_{2i}$, $i\in\N$, extends to an isometric isomorphism between $X$ and $Y:=\overline{\spa}\{e_{2i}: i\in\N\}$.
\item\label{main:3} The subspace $Y=\overline{\spa}\{e_{2i}: i\in\N\}$ is not complemented in $X$.   
\item\label{main:1} The sequence $(e_{2i-1})_{i=1}^\infty\subset X$ is unconditional.
\item\label{main:4} The space $X$ contains no complemented subspaces with an unconditional basis. 
\end{enumerate}
\end{theorem}
The natural setting for a space $X$ with the above properties is the class of Banach spaces with strictly controlled family of complemented subspaces, 
namely the spaces of Gowers-Maurey type. The Gowers-Maurey technique, started in \cite{gm} and developed in \cite{gm2}, see also \cite{maurey}, yields constructions solving a list of long-standing questions on Banach space structure, as unconditional basic sequence problem (solved in negative by the construction of the first known HI, i.e. hereditarily indecomposable Banach space \cite{gm}), Banach hyperplane problem, Schr\"oder-Bernstein problem for Banach spaces. Of particular interest for us is a general construction scheme of \cite{gm2} discussed briefly below, yielding spaces which provide answers to some of problems listed above and further questions on operator ideals and incomparability notions concerning Banach spaces, see   \cite{aiena:gonzales, ferenczi, horst:laustsen}. 
Introducing the Gowers-Maurey technique to the setting of Bourgain-Delbaen spaces by S.A.Argyros and R.Haydon \cite{argyros:haydon}, brought rapid development of theory of Banach spaces with small operator algebra, starting with the celebrated Argyros-Haydon space with "scalar+compact" property, towards a list of spaces with prescribed Calkin algebra, of which we mention here Tarbard's space with Calkin algebra isometric to the convolution algebra $\ell_1(\N_0)$ \cite{tarbard} that employs Gowers-Maurey scheme of \cite{gm2} mentioned above.  

We recall the general scheme of \cite{gm2}, as a starting point for the sequel. For two infinite $A=(a_i)_{i=1}^\infty,B=(b_i)_{i=1}^\infty\subset\N$,  a spread $S_{A,B}$ is a linear mapping $S_{A,B}:c_{00}\to c_{00}$ carrying each basic vector $e_{a_i}$ to $e_{b_i}$, $i\in\N$, and each basic vector $e_i$, $i\not\in A$, to 0. 
In \cite{gm2} for any so-called proper family of spreads $\mathcal{S}$ a Banach space $X(\mathcal{S})$ - a completion of $c_{00}$ under certain norm, is constructed so that each element  $S_{A,B}\in\mathcal{S}$ extends to a partial isometry on $X$ (carrying isometrically $\spa\{e_i:i\in A\}$ to $\spa\{e_i:i\in B\}$, with $\overline{\spa}\{e_i:i\in A\}$ and $\overline{\spa}\{e_i:i\in B\}$ 1-complemented in $X(\mathcal{S})$). What is crucial, Gowers-Maurey tools guarantee that, up to strictly singular perturbation, the algebra of bounded operators on $X(\mathcal{S})$ is generated by $\mathcal{S}$. If $\mathcal{S}$ is generated by $S_{\N,2\N}$, the resulting  space $X(\mathcal{S})$ is isometric to its subspace $Y=\overline{\spa}\{e_{2i}: i\in\N\}$ with $Y$ 1-complemented in $X(\mathcal{S})$. We  perturb the Gowers-Maurey construction, keeping the spread $S_{\N,2\N}$ as an isometry, but eliminating the continuity of the  adjoint spread $S_{2\N,\N}$ and projection $S_{2\N,2\N}$. As a by-product we obtain also an unconditional subsequence of the basis of the resulting space, contrary to the case of spaces $X(\mathcal{S})$ of \cite{gm2} containing no unconditional basic sequences. Instead of following directly  \cite{gm,gm2} we built the space $X$ within the classical by now framework of mixed Tsirelson spaces that was initiated in \cite{argyros:deliyanni} and discussed in detail in \cite{argyros:todorcevic}. 

The paper is organized as follows. In Section 2 we set the standard terminology, including mixed Tsirelson spaces and present the construction of the space $X$ with a basis $(e_i)_{i=1}^\infty$ (for the final step see Definition~\ref{def-X}), which yields immediately a 1-equivalence of
$(e_i)_{i=1}^\infty$ and $(e_{2i})_{i=1}^\infty$ (Proposition~\ref{isometric}), as well as unconditionality of $(e_{2i-1})_{i=1}^\infty$ (Proposition~\ref{prop:odd unc}). We recall the classical Gowers-Maurey tools following directly \cite{argyros:todorcevic} in Section 3.1.,  concluding it with the result on reflexivity of $X$ (Corollary \ref{reflexive}). Section 3.2. contains the analysis of the key tool - dependent sequences, applied in Section 3.3 for results on projections in $X$ 
(Theorem \ref{not-complemented} and Theorem \ref{unc-not-complemented}). The results listed above proves that the space $X$ defined in Section~\ref{section 2} satisfies assertion of  Theorem~\ref{main}.

\section{The construction of the space $X$}
\label{section 2}

We recall the standard notation.  
$\N$ denotes the set of positive integers, $\N_0:=\N\cup\{0\}$. For any  $A\subset\R$ and $k\in\N$ write $kA:=\{ka: a\in A\}\subset \R$ and $A+k:=\{a+k:a\in A\}\subset \R$.  $c_{00}$ denotes the vector space of sequences in $\R$ which are eventually  zero. We treat the elements of $c_{00}$ either as vectors, denoted by $x,y,\dots$, or as functionals on $c_{00}$, denoted by $f,g,\dots$. By $(e_i)_{i=1}^\infty$ we denote the unit vector basis of $c_{00}$, $(e_i^*)_{i=1}^\infty$ - the sequence of coordinate functionals. Given $v=(a_i)_{i=1}^\infty\in c_{00}$ we define
\begin{itemize}
    \item the support $\supp (v):=\{i\in\N: a_i\neq 0\}\subset\N$,
    \item the range $\range(v)\subset\N$ as the smallest interval in $\N$ containing $\supp(v)$,
    \item an element $|v|\in c_{00}$ with coordinates $|v|_i:=|v_i|$, $i\in\N$,
    \item the restriction $Ev\in c_{00}$ to $E\subset\N$ with coordinates $(Ev)_i=a_i$ if $i\in E$ and $(Ev)_i=0$ otherwise. 
\end{itemize} We write $k<v$ for $v\in c_{00}$, $k\in\N$, if either $v=0$ or $k<\min \supp (w)$ and $v<w$, for $v,w\in c_{00}$, if either $w=0$ or $v=0$ or $\max\supp v<\min \supp w$. A block sequence is any sequence $(v_n)_n\subset c_{00}$ with $v_{1}<v_{2}<\dots$. Notice that we allow in a block sequence null vectors - for the purpose of convenience in the notation in the sequel. 

\subsection{Mixed Tsirelson spaces} Spaces of Gowers-Maurey type are built on the basis of  mixed Tsirelson spaces, the theory of which started with the Schlumprecht space \cite{schlumprecht}, the first known arbitrary distortable space. We recall  the definition and basic properties of mixed Tsirelson spaces, defined by low complexity families $(\mathcal{A}_n)_{n=1}^\infty$, and spaces built on their basis, after \cite{argyros:todorcevic}. 

Fix two increasing sequences $(n_j)_{j=1}^\infty,(m_j)_{j=1}^\infty\subset\N$ with $m_j<n_j$, $j\in\N$. 

\begin{definition}[Mixed Tsirelson space]\label{def-mixed-tsirelson}  Let $\widetilde{K}\subset c_{00}$ be the smallest set satisfying the following:
\begin{enumerate}
    \item $(\pm e_i^*)_{i=1}^\infty\subset \widetilde{K}$,
    \item for any $j\in\N$ and $f_1<\dots<f_d$ in $\widetilde{K}$ with $d\leq n_j$ we have $m_j^{-1}(f_1+\dots+f_d)\in \widetilde{K}$.
\end{enumerate}
Let $\|\cdot\|$ be the norm on $c_{00}$ with the norming set $\widetilde{K}$, i.e. $\|\cdot\|=\sup\{|f(\cdot)|:f\in \widetilde{K}\}$. The mixed Tsirelson space
$T[(\mathcal{A}_{n_j},m_j^{-1})_{j\in\N}]$ is the completion of $(c_{00}, \|\cdot\|)$. 
\end{definition}

It is well-known that  the space $T[(\mathcal{A}_{n_j},m_j^{-1})_{j\in\N}]$ is  reflexive with $(e_i)_{i=1}^\infty$ as a normalized 1-unconditional basis   \cite[Thm I.10]{argyros:todorcevic}. The set $\widetilde{K}$ is symmetric (i.e. $f\in \widetilde K$ iff $-f\in\widetilde K$)  and closed under restrictions to arbitrary subsets (i.e. for any $f\in\widetilde K$ and $E\subset\N$ also $Ef\in \widetilde K$). 

\begin{remark}\label{cite-at}
For the future reference we note that in the terminology of \cite[Section II.1]{argyros:todorcevic}, the set $\widetilde K$ is equal to $W'_\kappa[G]$, with $\kappa=\omega$, a ground set $G=\{\pm e^*_i: i\in\N\}$, and $Y_G=c_0$. 
\end{remark}

Any $f\in\widetilde{K}$ of the form $f=m_j^{-1}(f_1+\dots+f_d)$ is called a \textbf{weighted functional}, with the \textbf{weight} $\w(f)=m_j^{-1}$. We note that the weight of a functional in $\widetilde{K}$ is not necessarily uniquely determined. We recall  the standard notion of a ree-analysis of a functional in $\widetilde K$. A tree  is a non-empty partially ordered  set $(\mt, \preceq)$ for which the set $\{ y \in \mt:y \preceq x \}$ is linearly ordered and finite for each $x \in \mt$. The {root} is the smallest element of the tree, terminal  nodes are the maximal elements. For any non-terminal $t\in\mt$,   $\suc(t)$  denotes the set of immediate successors  of $t$, i.e. all nodes $r\succ t$ with no $s \in \mt$ satisfying $r \succ s \succ t$.

\begin{definition}Let $f\in \widetilde K$. By a \textbf{tree-analysis} of
$f$ we mean a finite family $(f_t)_{t\in \mt}\subset \widetilde K$ indexed by a finite tree $\mt$ with a
unique root $\emptyset\in \mt$ such that 
\begin{enumerate}[(i)]
\item $f_\emptyset=f$,
\item $t\in \mt$ is terminal iff  $f_t\in \{\pm e_i^*:i\in\N\}$,
\item for every non-terminal $t\in \mt$,  $f_t=\pm m_j^{-1}E\sum_{r\in \suc(t)}f_r$, for some $j\in\N$ with $\# \suc(t)\leq n_j$, an interval $E\subset \N$, and a sequence $(f_r)_{r\in \suc(t)}\subset \widetilde K$ which is block with respect to an appropriate linear ordering of $\suc(t)$.
\end{enumerate}
\end{definition}
Note that by the construction any element of $\widetilde K$ admits a tree-analysis in $\widetilde K$.

\subsection{Mappings on $c_{00}$} We introduce here the mappings that will be used in the definition of a space, namely the spread $S$, its adjoint $R$.

Let $Y_0=\mathrm{span}\{e_{2i}:i\in\N\}\subset c_{00}$ and $S:c_{00}\to Y_0\subset  c_{00}$ be the spread of the form
\[  S: c_{00}\ni \sum_{i=1}^\infty a_ie_i\mapsto \sum_{i=1}^\infty a_ie_{2i}\in Y_0\subset c_{00}  \]

We shall use also the operator adjoint of $S$ restricted to $c_{00}$, $R: c_{00}\ni f\mapsto f\circ S\in c_{00}$, with  $c_{00}$ treated as a vector spaces of functionals, i.e. the linear map of the form
\[R: c_{00}\ni \sum_{i=1}^\infty b_ie_i^*\mapsto \sum_{i=1}^\infty b_{2i}e_i^*\in c_{00}\]
In the notation of \cite{gm2}, $R=S_{2\N,\N}$ and $S=S_{\N,2\N}$. 
We gather basic properties of the map $R$ below. 
\begin{fact}  \label{R properties} 
\begin{enumerate}[(1)]
    \item\label{R properties:ranges} For any $f\in c_{00}$ we have $R(2\range(Rf)f)=Rf$ and $2\supp(Rf)\subset\supp(f)$. 
    \item\label{R properties:ranges precise} For any  $g\in c_{00}$ with $\min\range(f),\max\range(f)\in 2\N$, $2\range(Rf)=\range(f)$.
       \item\label{R properties:S R blocks} $S$ and $R$ preserve block sequences, i.e. $Sx_1<Sx_2$ for any $x_1<x_2$ and $Rf_1<Rf_2$ for any $f_1<f_2$.
\end{enumerate}
\end{fact}

\subsection{The definition of the space $X$}

We shall use the Gowers-Maurey technique, introduced in \cite{gm} and extended in \cite{gm2} within the framework of \cite{argyros:todorcevic}. 
\begin{assumption}\label{assumption-integers}
Fix two sequences $(n_j)_j, (m_j)_j\subset 2\N$ satisfying the following: $m_1=2$, $n_1= 4$, $m_{j+1}=m_j^5$, $n_{j+1}= (5n_j)^{s_j}$, where $s_j=\log_2(m_{j+1}^3)$, for all $j\in\N$. 
\end{assumption}
Denote by $c_{00}^\mathbb{Q}$ the $\mathbb{Q}$-vector subspace of $c_{00}$ spanned by $\{e_i^*: i\in\N\}$ and fix  an injective function $\sigma: \{(f_1,\dots,f_d)\subset c_{00}^\mathbb{Q}: f_1<\dots<f_d\}\to 4\N$ such that 
\[\sigma((f_1,\dots,f_d))\geq 4\max\supp (f_d)\|f_1+\dots+f_d\|_\infty^{-1},\quad (f_1,\dots,f_d)\subset c_{00}^\Q.\]

A block sequence $f_1<\dots<f_{2d}$ in $\widetilde K$, $d\leq n_{2j-1}/2$, $j\in \N$, is called \textbf{$j$-special} with associated sequence $(k_1,\dots,k_d)\subset\N$, if
\begin{enumerate}[(S1)]
    \item\label{s-1} $\w(f_1)=m_{4l-2}^{-1}$ for some $l\in\N$ with $m_{4l-2}>9n_{2j-1}^2$,
    \item\label{s-2} $\w(f_{2i+2})=\w(f_{2i+1})=m^{-1}_{\sigma((f_1,f_2,\dots,f_{2i-1},f_{2i}))}$ for any $i=1,\dots, d-1$,
    \item\label{s-3} $f_{2i-1}=R^{k_i}f_{2i}$ with $\range (f_{2i})=2^{k_i}\range(f_{2i-1})$, for any $i=1,\dots,d$.
\end{enumerate} 

A block sequence $g_1<\dots<g_{2d}$ in $\widetilde K$, with $d\leq n_{2j-1}/2$, $j\in\N$,  $m\in\Z$, is called a \textbf{shifted $j$-special sequence $m$-modeled on a $j$-special sequence} $f_1<\dots<f_{2d}$ with associated sequence $k_1,\dots,k_d\in\N$,  provided
\begin{enumerate}[(SS1)]
    \item\label{ss-1} $\w(g_i)=\w(f_i)$, for all $i=1,\dots,2d$ with $g_i\neq 0$,
    \item\label{ss-2} $\range (g_i)\subset 2^m\range(f_i)$, for all $i=1,\dots,2d$, 
    \item\label{ss-3}  $g_{2i-1}=R^{k_i}g_{2i}$, for all $i=1,\dots,d$.
\end{enumerate}
By definition, $j$-special sequences are shifted $j$-special sequences. 
\begin{remark}\label{special sequences properties}
\begin{enumerate}[(1)]
\item\label{special sequences properties:0} Note that we allow in a shifted special sequence null vectors. However, for a shifted  special sequence $g_1<\dots<g_{2d}$, if $g_{2i-1}\neq 0$ then also $g_{2i}\neq 0$, $i=1,\dots,d$.
\item\label{special sequences properties:2} The weight of a functional in $\widetilde K$ is not uniquely defined. However, for elements of a shifted special sequence, weights will refer to those  in \ref{s-1}, \ref{s-2} or \ref{ss-1}.
\item\label{special sequences properties:3}  Note the tree-like property of the family of special sequences: for any special sequences $f_1<\dots<f_{2d}$ and $h_1<\dots<h_{2a}$ and $r:=\min\{i=1,\dots,\min\{d,a\}, f_{2i-1}\neq h_{2i-1}\}$ we have $\w(f_i)\neq \w(h_s)$ for any $(s,i)\in\{2r+1,\dots,2a\}\times\{1,\dots,2d\}$. 
\item\label{special sequences properties:4} 
Fix a shifted $j$-special sequence  $(g_1,\dots,g_{2d})$, $j\in\N$. Then  $(Rg_1,\dots,Rg_{2d})$ is a shifted $j$-special sequence 
(see Fact~\ref{R properties}\ref{R properties:ranges} and \ref{R properties:ranges precise}).
\end{enumerate}
\end{remark}

\noindent\textbf{Definition of the norming set $K$ for the space $X$.} The norming set $K$ of the space $X$ is defined as the union of an increasing sequence of sets $K_n\subset \widetilde K$, $n\in\N_0$, with $\widetilde K$ introduced in Definition~\ref{def-mixed-tsirelson}, which we define inductively below, verifying at each step the following properties.
\begin{enumerate}[(K1)]
\item\label{k:1} $K_n\subset\widetilde K$, $K_n$ is symmetric and closed under restriction to intervals. 
\item\label{k:2} Any $f\in K_n$ has a tree-analysis $(f_t)_{t\in\mt}\subset K_n$. 
\item\label{k:3} $R(K_n)=K_n$. Moreover, for any weighted $f\in K_n$ there is $g\in K_n$ with $Rg=f$ and $\w(g)=\w(f)$. 
\end{enumerate}
Let $K_0=\{\pm e^*_i: i\in \N\}\cup\{0\}$ and note that $K_0$ satisfies \ref{k:1}-\ref{k:3}. Fix $n\in\N$ and assume that  $K_{n-1}$ is already defined and satisfies \ref{k:1}-\ref{k:3}.  
Let $K_n$ be the collection of all functionals $f$ of one of the forms listed below. 
\begin{enumerate}[(1)]
\item\label{inductive-K:1} $f\in K_{n-1}$. 
\item\label{inductive-K:2} $f=m_{2j}^{-1}(f_1+\dots+f_d)$, with $j\in\N$ and a block sequence $f_1<\dots<f_d$, $d\leq n_{2j}$, of elements of  $K_{n-1}$. Such $f$ is called a \textbf{regular} functional. 
\item\label{inductive-K:3} $f=m_{2j-1}^{-1}E(g_1+\dots+g_{2d})$, for $j\in\N$, an interval $E\subset\N$ and a shifted $j$-special sequence  $g_1<\dots<g_{2d}$, $d\leqslant n_{2j-1}/2$, of elements of $K_{n-1}$ modeled on a $j$-special sequence of elements on $K_{n-1}$. Such $f$ is called a \textbf{special} functional. 
\end{enumerate}
Note that $K_n\subset \widetilde K$ by linearity of $R$ and the inductive assumption. We verify conditions \ref{k:1}-\ref{k:3} for $K_n$. 

\ref{k:1}. The operation "restriction to an interval" in \ref{inductive-K:2} is not included, contrary to the case of \ref{inductive-K:3}, as by the inductive assumption, for any regular $f$ and an interval $E\subset\N$,  $Ef$ is regular as well. Therefore \ref{k:1} follows by the definition of $K_n$ and the inductive assumption.  

\ref{k:2} follows immediately by the definition of $K_n$ and properties \ref{k:2} and \ref{k:3} of $K_{n-1}$. 

\ref{k:3}: the inclusion $R(K_n)\subset K_n$. For functionals from $K_{n-1}$ or regular functionals use \ref{k:1} and \ref{k:3} for $K_{n-1}$, for special functionals use Remark~\ref{special sequences properties}\ref{special sequences properties:4} and \ref{k:3} for $K_{n-1}$. 

\ref{k:3}: the inclusion $K_n\subset R(K_n)$, and the "moreover" part. For $f\in K_{n-1}$ use \ref{k:3} for $K_{n-1}$. Assume $f\in K_n$ is regular, i.e. $f=m_{2j}^{-1}(f_1+\dots+f_d)$, with $f_1<\dots<f_d$ in $K_{n-1}$. By \ref{k:3} and \ref{k:1} for $K_{n-1}$ pick $h_i\in K_{n-1}$ with $Rh_i=f_i$ and $\range h_i=2\range f_i$, $i=1,\dots,d$. Then $h:=m_{2j}^{-1}(h_1+\dots+h_d)\in K_n$ satisfies $Rh=f$.

Assume $f\in K_n$ is special, i.e. of the form $f=m_{2j-1}^{-1}E(g_1+\dots+g_{2d})$, for $j\in\N$, an interval $E\subset\N$ and a shifted $j$-special sequence  $g_1<\dots<g_{2d}\subset K_{n-1}$, $d\leqslant n_{2j-1}/2$,  $m$-modeled on a $j$-special sequence $f_1<\dots<f_{2d}\subset K_{n-1}$, with associated $k_1,\dots,k_d\subset\N$ as in \ref{s-3}. For a fixed $i=1,\dots,d$  choose $h_{2i-1}$ and $h_{2i}$ as follows. Let $h_{2i-1}:=R^{k_i-1}g_{2i}\in K_{n-1}$ by \ref{k:3} for $K_{n-1}$, $i=1,\dots,d$. Then $Rh_{2i-1}=g_{2i-1}$ and, if $h_{2i-1}\neq 0$, $\w (h_{2i-1})=\w(g_{2i})$. Moreover, $2^{k_i-1}\range (h_{2i-1})\subset \range(g_{2i})\subset 2^m\range(f_{2i})=2^{m+k_i}\range(f_{2i-1})$ by \ref{s-3}. If $g_{2i}\neq 0$ pick $h_{2i}\in K_{n-1}$ with $Rh_{2i}=g_{2i}$, $\w(h_{2i})=\w(g_{2i})$ and $\range(h_{2i})=2\range(g_{2i})$ by \ref{k:3} and \ref{k:1} for $K_{n-1}$. It follows that $h_{2i-1}=R^{k_i}(h_{2i})$ and $\range(h_{2i})\subset 2^{m+1}\range(f_{2i})$. Therefore, $h_1<\dots<h_{2d}$ is a shifted $j$-special sequence of elements of $K_{n-1}$, $(m+1)$-modeled on $f_1<\dots<f_{2d}$, thus for a special functional $h:=m_{2j-1}^{-1}(2E)(h_1+\dots+h_{2d})\in K_n$ we have  $f=Rh$, which ends the proof.

\begin{definition}[The space $X$]\label{def-X}
Let $K=\bigcup_{n=0}^\infty K_n\subset\widetilde K$, with  $\widetilde K$ introduced in Definition~\ref{def-mixed-tsirelson} and $(K_n)_{n=0}^\infty$ defined above. Let $\|\cdot\|$ be the norm on $c_{00}$ defined by $K$ as the norming set, i.e. $\|\cdot\|=\sup\{|f(\cdot)|:f\in K\}$. Finally, let $X$ be the completion of $(c_{00}, \|\cdot\|)$ and let $Y=\overline{Y_0}\subset X$. 
\end{definition}
\begin{notation}\label{notation-X} We extend to the setting of $X$ the notions of a support and  a range of a vector, therefore also of a block sequence. We call a subspace $W$ of $X$ a block subspace, provided $W$ is a closed vector subspace spanned by a block sequence. Finally, we denote by $X_m$, $m\in\N$, the "tail" subspace of $X$; that is $X_m=\{x\in X: x> m\}$, $m\in\N$.
\end{notation}
\begin{fact}\label{fact after def X}
 \begin{enumerate}[(1)]
\item\label{fact after def X:1} As $K_0\subset K$, $\|\cdot\|_\infty\leq \|\cdot\|$. 
\item\label{fact after def X:2} $K$ satisfies properties \ref{k:1}-\ref{k:3}.
\item\label{fact after def X:3} By \ref{k:1}, the unit vectors $(e_i)_{i=1}^\infty$ form a bimonotone basis of $X$. 
\item\label{fact after def X:4} $m_{2j}^{-1}(f_1+\dots+f_d)\in K$ for any elements $f_1<\dots<f_d$ of $K$ with  $d\leq n_{2j}$, $j\in\N$. 
 \end{enumerate}   
\end{fact}

\begin{remark}\label{remark after def X}
In the construction of \cite{gm2} the norming set was closed both under $R$ and $S$, which gives a bounded projection on $Y$ in the resulting space. In order to avoid it, we close the norming set under $R$ and some selections of  $R^{-1}$, typically different from $S$ (cf. conditions \ref{s-3} and \ref{ss-3} in the definition of special sequences and shifted special sequences). This precise choice of selections of $R^{-1}$ admitted in $K$, by the Gowers-Maurey technique, yields uncomplementability of  $Y$ and of subspaces with unconditional basis in $X$ (see Section \ref{section 3}). At the same time the flexibility of the definition of shifted special sequences allows for unconditionality of the subsequence $(e_{2i-1})_{i=1}^\infty$ of the basis of $X$. 
\end{remark}

\begin{proposition}\label{isometric}
The injection $S:c_{00}\to Y_0\subset c_{00}$ is an isometry with respect to the norm $\|\cdot\|$, thus extends to an isometry $X\to X$ with $S(X)=Y$. In other words, $(e_i)_{i=1}^\infty\subset X$ and $(e_{2i-1})_{i=1}^\infty\subset Y$ are 1-equivalent.
\end{proposition}
\begin{proof}
Fix $x\in c_{00}$. 
Given $f\in K$, by \ref{k:3} pick $g\in  K$ with $f=Rg$ and estimate $|f(x)|=|(Rg)(x)|=|g(Sx)|\leq \|Sx\|$. On the other hand, for any $g\in K$ also $Rg\in K$ by \ref{k:3}, thus $|g(Sx)|=|(Rg)(x)|\leq \|x\|$. Thus  $\|x\|=\|Sx\|$ for any $x\in c_{00}$, which ends the proof, as $Se_i=e_{2i}$, $i\in\N$, by the very definition. 
\end{proof}

Recall that a basis  $(v_n)_{n=1}^\infty$ of a Banach space $V$ is $C$-\textbf{unconditional}, $C\geqslant 1$, if for any signs $\xi=(\varepsilon_i)_{i=1}^\infty\subset\{\pm 1\}$, the map carrying each $v_i$ to $\varepsilon_iv_i$, $i\in\N$, extends to a bounded  linear isomorphism of $V$ of norm at most $C$. A basis $(v_n)_{n=1}^\infty$ of a Banach space is unconditional, if it is $C$-unconditional for some $C\geqslant 1$.  We shall prove  that the sequence $(e_{2i-1})_{i\in\N}\subset X$ is unconditional, based on Lemma~\ref{lemma:odd unc}. We use a standard notation: for  $f\in K$ and $A\subset\N$ write $f\cdot\chi_A$ for a functional on $c_{00}$ with $(f\cdot\chi_A)(e_i)=f(e_i)$ for $i\in A$ and $(f\cdot\chi_A)(e_i)=0$ otherwise.
\begin{lemma}\label{lemma:odd unc}
Take $f\in K$ and fix $n\in 2\N-1$. Let $g=f\cdot\chi_J$, where $J=\N\setminus\{2^kn: k\in \N\cup\{0\}\}$. Then $g\in K$ with $g\cdot\chi_{2\N-1}=f\cdot\chi_{(2\N-1)\setminus\{n\}}$.
\end{lemma}
\begin{proof} The second part is clear by the definition of $g$, for the first part let $(f_t)_{t\in\mt}\subset K$ be a tree-analysis of $f$. We prove inductively on $t\in\mt$, starting from terminal nodes, that $g_t=f_t\cdot\chi_J\in K$, $t\in\mt$. In particular, $g=g_\emptyset\in K$, which ends the proof. 

For terminal $t\in\mt$ the claim is clear. Take non-terminal $t\in\mt$. If $f_t$ is regular, i.e. of the form $f_t=m_{2j}^{-1}\sum_{s\in\suc(t)}f_s$, $j\in\N$, then by the inductive assumption $g_s=f_s\cdot\chi_J\in K$, $s\in\suc(t)$, thus $g_t=m_{2j}^{-1}\sum_{s\in\suc(t)}g_s\in K$ by Fact~\ref{fact after def X}\ref{fact after def X:4}. Assume $f_t$ is special, i.e. $f_t=Em_{2j-1}^{-1}\sum_{i=1}^{2d}\widetilde f_i$, with $(f_s)_{s\in \suc(t)}=(\widetilde f_i)_{i=1}^{2d}$, $d\leqslant n_{2j-1}/2$, $j\in\N$, for some shifted $j$-special sequence $\widetilde f_1<\dots<\widetilde f_{2d}$ with associated $k_1,\dots,k_d\in\N$. Enumerate accordingly $(g_s)_{s\in\suc(t)}=(\widetilde g_i)_{i=1}^{2d}$. By the inductive assumption and the definition of $(g_s)_{s\in\suc(t)}$, $\widetilde g_i\in K$, $\supp \widetilde g_i\subset\supp\widetilde f_i$, $\w(\widetilde g_i)=\w(\widetilde f_i)$, $i=1,\dots,2d$. Moreover, as $\widetilde f_{2i-1}=R^{k_i}\widetilde f_{2i}$, also $ \widetilde g_{2i-1}=R^{k_i}\widetilde g_{2i}$, $i=1,\dots,d$  (by the choice of $n\in 2\N-1$ and $J$). Therefore $\widetilde g_1<\dots<\widetilde g_{2d}$ is a shifted $j$-special sequence and thus $g_t=Em_{2j}^{-1}\sum_{s\in\suc(t)}g_s\in K$, which concludes the inductive step. 
\end{proof}

\begin{proposition}\label{prop:odd unc}
    The subsequence $(e_{2i-1})_{i\in\N}$ of the basis of $X$ is unconditional. 
\end{proposition}
\begin{proof} Take any $f\in K$ and $I\subset 2\N-1$. As $f\in c_{00}$, iterating the procedure of Lemma~\ref{lemma:odd unc} we find $g\in K$ with $g\cdot\chi_{2\N-1}=f\cdot\chi_I$. Consequently, with $Z:=\overline{\spa}\{e_{2i-1}:i\in\N\}\subset X$, any coordinate projection 
$P_I:Z\ni \sum_{n=1}^\infty a_ne_{2n-1}\to \sum_{n\in I}a_ne_{2n-1}\in Z$, $\emptyset\neq I\subset 2\N-1$, is well-defined and  bounded with $\|P_I\|=1$, which ends the proof.
\end{proof}

\section{Projections in the space $X$}
\label{section 3}
\subsection{Classical tools in $X$}

We recall here the notions and results from \cite{argyros:todorcevic}. First note, that $K$ ($=D_G$ in the notation of \cite{argyros:todorcevic}) by properties \ref{k:1}, \ref{k:2} and  Fact \ref{fact after def X}\ref{fact after def X:4}, is a tree-complete extension of the ground set $G:=\{\pm e^*_i: i\in\N\}$ (see \cite[Definition~II.10]{argyros:todorcevic}), therefore we can directly cite results from \cite{argyros:todorcevic}. The space $X$ defined in the previous section is equal to the space $\mathfrak{X}_\omega(D_G)$  of \cite[Section II.2]{argyros:todorcevic}. Recall that a bounded operator $T$ between Banach spaces $U,V$ is \textbf{strictly singular}, if any infinite dimensional subspace $W$ of $U$ contains, for any $\varepsilon>0$, a vector $w$ with $\|w\|>1$ and $\|Tw\|<\varepsilon$. 
\begin{proposition}\label{cite-at-ss-extension}
The formal identity map $X\to c_0$ is strictly singular. 
\end{proposition}
\begin{proof}
By the Krein--Milman--Rutman Theorem, when checking strict singularity of the operator on $X$ we can consider only block subspaces in $X$. In a given block subspace, for $j\in\N$, pick a normalized block sequence $x_1<\dots<x_{n_{2j}}$ and a block sequence  $f_1<\dots<f_{n_{2j}}$ in $K$ with $f_i(x_i)=1$ for any $i=1,\dots, n_{2j}$. Let $x=m_{2j}n_{2j}^{-1}(x_1+\dots+x_{n_{2j}})$. Then, by the definition of $K$, also $f=m_{2j}^{-1}(f_1+\dots+f_{n_{2j}})\in K$, thus $\|x\|\geq f(x)=1$, whereas \[\|x\|_\infty=m_{2j}n_{2j}^{-1}\max\{\|x_i\|_\infty: i=1,\dots,n_{2j}\}\leq m_{2j}n_{2j}^{-1}\max\{\|x_i\|, i=1,\dots,n_{2j}\}= m_{2j}n_{2j}^{-1}\] (see Fact \ref{fact after def X}\ref{fact after def X:1}). By Assumption \ref{assumption-integers},  $m_{2j}n_{2j}^{-1}\to 0$, as $j\to\infty$, which ends the proof.      
\end{proof}
The above observation stated in the terminology of \cite[Section II.2]{argyros:todorcevic} claims that $X$ is a strictly singular extension of $Y_G(=c_0)$ (cf. Remark \ref{cite-at}).

\begin{definition}\cite[Definition~II.21]{argyros:todorcevic}
A normalized vector $x\in X$  is called a $C$-$\ell_1^N$-\textbf{average}, $C>1$, $N\in\N$, if $x=(x_1+\dots+x_N)/N$ for a block sequence $(x_n)_{n=1}^N\subset X$ with $\max_n\|x_n\|\leq C$.  
\end{definition} 
If $C=2$ we  omit it  and say  that $x$ is an $\ell_1^N$-\textbf{average}. 
By Proposition~\ref{cite-at-ss-extension} we have the following 
\begin{lemma}\label{av-existence}\cite[Lemma~II.22]{argyros:todorcevic}  Any block subspace of $X$ for any $N\in\N$ and $\varepsilon>0$ contains an $\ell_1^N$-average $x$ with $\|x\|_\infty<\varepsilon$. 
\end{lemma}

\begin{lemma}\label{av-est}\cite[Lemma~II.23]{argyros:todorcevic}
Let $x\in X$ be an $\ell_1^{n_j}$-average, for some $j\in\N$. Then for any $f\in K$ with $\w(f)>m_j^{-1}$ we have
$|f(x)|\leq 3\w(f)$.
\end{lemma}

\begin{definition}[RIS]\cite[Definition~II.13]{argyros:todorcevic}\label{def-ris} A block sequence $(x_i)_{i\in I}\subset
X$ is called a $(C,\varepsilon)$- \textbf{rapidly increasing sequence} (RIS), for $C>1$, $\varepsilon>0$ and an interval $I\subset\N$, if there is a strictly increasing
sequence $(j_i)_{i\in I}\subset L_0$ such that 

\begin{enumerate}[(1)]

\item\label{def-ris:1} $\|x_i\|\leq C$, $\|x_i\|_\infty<\varepsilon$, for any $i\in I$,

\item\label{def-ris:2} 
$2m_{j_{i+1}}^{-1}\#\supp(x_i)<\varepsilon$ for any $i\in I\setminus\{\max I\}$,

\item\label{def-ris:3} $| f(x_i)|\leq C\w(f)$ for any weighted $f\in K$ with $\w(f)> m_{j_i}^{-1}$ for any $i\in I$.

\end{enumerate}
\end{definition}

Following the proof of \cite[Proposition~II.25]{argyros:todorcevic} we obtain the following.
\begin{lemma}\label{ris-existence}
Any block subspace of $X$ contains for any $\varepsilon>0$ a normalized block sequence $(x_i)_{i=1}^\infty$ 
such that   $(S^kx_i)_{i=1}^\infty$ is a $(3,\varepsilon)$-RIS for any $k\in\N_0$.
\end{lemma}
\begin{proof} 
In a fixed block subspace of $X$  pick inductively (using Lemma~\ref{av-existence}) a block sequence $(x_i)_{i=1}^\infty$ of $\ell_1$-averages and  $(j_i)_{i=1}^\infty\subset \N$ satisfying  \ref{def-ris:2} of Definition~\ref{def-ris} such that $x_i$ is an $\ell_1^{j_i}$-average with $\|x_i\|_\infty<\varepsilon$, for all $i\in\N$. Note that for any $k\in\N$,  $(S^kx_i)_{i=1}^\infty$ is also a block sequence of $\ell_1^{j_i}$-averages (as $S$ is an isometry), that satisfies \ref{def-ris:1} and \ref{def-ris:2} of Definition~\ref{def-ris}. By Lemma~\ref{av-est} each sequence $(S^kx_i)_{i=1}^\infty$, $k\in\N_0$, satisfies  \ref{def-ris:3} of Definition~\ref{def-ris}, thus is a $(3,\varepsilon)$-RIS. 
\end{proof}

\begin{proposition}\cite[Proposition~II.19]{argyros:todorcevic}\label{basic-inequality-cor} Fix $j\in\N$. Let $(x_i)_{i=1}^d\subset X$, $d\leq n_j$, be a $(C,m_j^{-2})$-RIS. Then for any weighted functional $f\in K$  we have
$$
\Big| f\Big(\frac{1}{n_j}\sum_{i=1}^dx_i\Big)\Big|\leq \begin{cases} \frac{2C}{m_j^2} \ \ \ \ \ \ \ \ \text{if}\ \ \ \w(f)<\frac{1}{m_j} \\
\frac{2C}{m_j} \ \ \ \ \ \ \ \ \text{if}\ \ \ \w(f)=\frac{1}{m_j}\\
\frac{3C}{m_j}\w(f) \ \ \ \text{if}\ \ \ \w(f)>\frac{1}{m_j}
\end{cases}
$$
In particular
$$
\Big\|\frac{1}{n_j}\sum_{i=1}^dx_i\Big\|\leq \frac{2C}{m_j}
$$
Moreover, if for any weighted  $h\in K$ with $\w(h)=m_j^{-1}$ and any interval $F\subset \{1,\dots,d\}$,
\[\Big|h\Big(\sum_{i\in F}x_i\Big)\Big|\leq C, \]
then
\[\Big\|\frac{1}{n_j}\sum_{i=1}^dx_i\Big\|\leq \frac{4C}{m_j^2}\]
\end{proposition}
\begin{corollary}
\label{reflexive}\cite[Proposition~II.28]{argyros:todorcevic} The basis $(e_i)_{i=1}^\infty$ of $X$ is boundedly complete and shrinking, thus the space $X$ is reflexive.   
\end{corollary}

\subsection{Dependent sequences in $X$} This section is devoted to the key tool in the analysis of bounded operators on $X$, namely dependent sequences. We start with the notion of an exact pair, the basic "brick" for building dependent sequences. 
\begin{definition}\cite[Definition~II.31]{argyros:todorcevic}\label{def-exact}
A pair $(x,f)\in X\times K$ is called a $j$-\textbf{exact pair}, $j\in \N$, if 
\begin{enumerate}[(1)]
\item\label{def-exact:1} $f$ is weighted with $\w(f)=m_j^{-1}$,
\item\label{def-exact:2} $1\leq \|x\|\leq 6$, $\|x\|_\infty<m_{2j}^{-1}$,
\item\label{def-exact:3} $f(x)=1$,  $\range(f)=\range(x)$,
\item\label{def-exact:4} $|h(x)|\leq 9\max\{m_{2j}^{-1},\w(h)\}$ for any weighted $h\in K$ with weight $\w(h)\neq m_{2j}^{-1}$. 
\end{enumerate}
\end{definition}
Repeating the proof of \cite[Proposition~II.32]{argyros:todorcevic} we obtain the following. 
\begin{lemma}\label{exact existence}
Any block subspace of $X$, for any $j\in \N$ and $\varepsilon>0$, contains a vector $x$ such that for some $f\in K$, $(x,f)$ is a  $j$-exact pair.
\end{lemma}
\begin{proof} In a given block subspace of $X$ pick (by Lemma~\ref{ris-existence}) a normalized block sequence $(x_i)_{i=1}^{n_{2j}}$ with   $(S^kx_i)_{i=1}^{n_{2j}}$ a $(3,m_{2j}^{-1})$-RIS for any $k\in\N_0$. For each $i=1,\dots,n_{2j}$ pick $f_i\in K$ with $f_i(x_i)=1$ and $\range(f_i)=\range (x_i)$ (by bimonotonicity of the basis in $X$). Let 
\[x=\frac{m_{2j}}{n_{2j}}\sum_{i=1}^{n_{2j}}x_i, \ \ \ f=\frac1{m_{2j}}\sum_{i=1}^{n_{2j}}f_i\in K\]
By Proposition~\ref{basic-inequality-cor}, $x$ satisfies \ref{def-exact:2} and \ref{def-exact:4}, which ends the proof. \end{proof}

\begin{definition}[Dependent sequences]\label{def dep} A double sequence $(x_i,f_i)_{i=1}^{n_{2j-1}}\subset
X\times K$ is called a $j$-\textbf{dependent sequence}, $j\in\N$, if for some  $(k_i)_{i=1}^{n_{2j-1}/2}, (j_i)_{i=1}^{n_{2j-1}/2}\subset 2\N$  
we have the following
\begin{enumerate}[(D1)]
\item\label{dep-seq-initial} 
$m_{j_1}>9n_{2j-1}^2$,
\item\label{dep-seq-special} $(f_i)_{i=1}^{n_{2j-1}}$ is a $j$-special sequence with associated $(k_i)_{i=1}^{n_{2j-1}/2}$,
 \item\label{dep-seq-exact} $(x_{2i-1},f_{2i-1})$ and $(x_{2i}, f_{2i})$ are $j_i$-exact pairs, for any $i=1,\dots,n_{2j-1}/2$,
\item\label{dep-seq-vectors-shifted} $x_{2i}=S^{k_i}(x_{2i-1})$,  $i=1,\dots,n_{2j-1}/2$,

\item\label{dep-seq-limit-intersection} for $k\in\N_0$ there is at most one $i\in\{1,\dots,n_{2j-1}/2\}$ with $2^k\range(x_{2i-1})\cap \range(x_{2i})\neq \emptyset$. 
\end{enumerate}
\end{definition}

\begin{remark}\label{remark dep}
Note that by \ref{dep-seq-initial}, \ref{dep-seq-special}, \ref{dep-seq-exact} and the choice of $\sigma$, we have \begin{enumerate}[(D6)]
        \item\label{dep-seq-estimate} $4m_{j_{i+1}}^{-1}\#\supp(x_{2i-1})<m_{2j-1}^{-2}$   for any $i=1, \dots, n_{2j-1}-1$. 
    \end{enumerate}

\end{remark}
\begin{proposition}\label{prop: choice of dep seq}
    For any block subspace $W\subset X$ and $j\in\N$ we have the following:
    \begin{align*}
        \forall & N_1\in\N\quad \exists (x_1,f_1)\in (W\cap X_{N_1})\times K:\\ 
        \forall & N_2\in\N\quad\exists (x_2,f_2)\in X_{N_2}\times K:  \\  &  \dots \\
        \forall & N_{n_{2j-1}-1}\in\N\quad \exists (x_{n_{2j-1}-1},f_{n_{2j-1}-1})\in (W\cap X_{N_{n_{2j-1}-1}})\times K:\\ 
        \forall & N_{n_{2j-1}}\in\N\quad\exists (x_{n_{2j-1}},f_{n_{2j-1}})\in X_{N_{n_{2j-1}}}\times K,
    \end{align*}    
    such that $(x_i,f_i)_{i=1}^{n_{2j-1}}$ is a $j$-dependent sequence. 
\end{proposition}
\begin{proof} 
    For a fixed $j\in\N$ and consecutive choice of $(N_i)_{i=1}^{n_{2j-1}}$ we shall pick inductively a double sequence $(x_i,f_i)_{i=1}^{n_{2j-1}}$, $(k_i)_{i=1}^{n_{2j-1}/2}, (s_i)_{i=1}^{n_{2j-1}/2}\subset\N$ satisfying \ref{dep-seq-initial}-\ref{dep-seq-vectors-shifted} and 
\begin{enumerate}[(D7)]
    \item\label{dep-seq-lemma}   $\range(x_{2i-1})\subset [2, 2^{s_i}]$ and $k_i+s_i<k_{i+1}-s_{i+1}$ for any $i=1,\dots,n_{2j-1}$. 
\end{enumerate}
    Given $N_1\in\N$, by Lemma~\ref{exact existence} pick a $j_1$-exact pair $(x_1,f_1)$ with $x_1\in W$, $x_1>N_1$ and $m_{j_1}>9n^2_{2j-1}$, pick $s_1\in\N$ to satisfy the first part of \ref{dep-seq-lemma} and, given $N_2\in\N$, choose $k_1\in\N$ so that $x_2=S^{k_1}(x_1)$ satisfies $x_2>x_1$ and $x_2>N_2$. By \ref{k:3} pick $f_2\in K$ with $\w(f_2)=\w(f_1)$, $R^{k_1}f_2=f_1$ and $\range (f_2)=2^{k_1}\range(f_1)$. Fix $i<n_{2j-1}/2$ and assume that $2i$ elements with desired properties are chosen. Given $N_{2i+1}\in\N$, for $j_{i+1}=\sigma(f_1,\dots,f_{2i})$ pick by Lemma~\ref{exact existence} a $j_{i+1}$-exact pair $(x_{2i+1},f_{2i+1})$ with $x_{2i+1}\in W$,   $x_{2i+1}>x_{2i}$ and $x_{2i+1}>N_{2i+1}$. Given $N_{2i+2}$ choose $s_{i+1},k_{i+1}>N_{2i+2}\in\N$ fulfilling \ref{dep-seq-lemma} and let $x_{2i+2}=S^{k_{i+1}}(x_{2i+1})$. By the property \ref{k:3} pick $f_{2i+2}\in K$ with $\w(f_{2i+2})=\w(f_{2i+1})$, $R^{k_i}f_{2i+2}=f_{2i+1}$ and $\range(f_{2i+2})=2^{k_i}\range(f_{2i+1})$, which completes the inductive step. 

To conclude the proof we verify  \ref{dep-seq-limit-intersection}. 
Fix $k\in\N$. For any  $i=1,\dots,n_{2j-1}$,  $\range(x_{2i})\subset [2^{k_i}, 2^{k_i+s_i}]$ and   $2^k\range(x_{2i-1})\subset [2^{k},2^{k+s_i}]$. Hence, if $2^k\range(x_{2i-1})\cap\range(x_{2i})\neq \emptyset$ for some $i=1,\dots,n_{2j-1}/2$, then $k\leq k_i+s_i$ and $k_i\leq k+s_i$, i.e. $k\in [k_i-s_i, k_i+s_i]$, and by \ref{dep-seq-lemma} there is at most one $i$ with such property, which ends the proof of \ref{dep-seq-limit-intersection}.  
\end{proof}
\begin{proposition}\label{dep-est} Let $(x_i,f_i)_{i=1}^{n_{2j-1}}$ be a $j$-dependent
sequence, $j\in\N$. Then
\[\Big\|\frac{1}{n_{2j-1}}\sum_{i=1}^{n_{2j-1}}(-1)^ix_i\Big\|\leq \frac{240}{m_{2j-1}^2}\]
\end{proposition}

\begin{proof}
Notice that the sequence $(x_{2i}-x_{2i-1})_{i=1}^{n_{2j-1}/2}$ is itself a $(18,m_{2j-1}^{-2})$-RIS by \ref{dep-seq-exact} and \ref{dep-seq-estimate}.

By Proposition~ \ref{basic-inequality-cor} it is enough to show that for any special $f\in K$ with $\w(f)=m_{2j-1}^{-1}$ and any interval $E$ we have $| f(\sum_{i\in F}(x_{2i}-x_{2i-1}))|\leq 60$. 

Fix an interval $F\subset\N$ and  $f\in K$ of the form $f=m_{2j-1}^{-1}E(g_{2l-1}+\dots+g_{2d})$, 
for an interval $E\subset\N$ and a shifted $j$-special sequence
$g_1<\dots<g_{2d}$ $m$-modeled on a $j$-special sequence $h_1<\dots<h_{2d}$ with $m\in\Z$ and $1\leq l< r<
d\leq n_{2j-1}/2$ and $E(g_{2d-1}+g_{2d})\neq 0$, where $l=\min\{i=1,\dots,d:E(g_{2i-1}+g_{2i})\neq 0\}$ and $r=\min\{i=1,\dots, d: h_{2r+1}\neq f_{2r+1}\}$ (for $l,r$ not well-defined or other relations between $l,r,d$ the estimates below simplify). 

Note the following estimates. 
For any $(s,i)\in (\{2r+3,\dots,d\}\times \{1,\dots,n_{2j-1}\})\cup (\{2r+1,2r+2\}\times \{2r+3,\dots,n_{2j-1}\})$, by Remark \ref{special sequences properties}\ref{special sequences properties:3},  $\w(g_s)=\w(h_s)\neq \w(f_i)$, thus $|
(g_s)(x_i)|\leq n_{2j-1}^{-2}$ by \ref{dep-seq-exact} and \ref{dep-seq-estimate}. For any $i\neq s$ by \ref{dep-seq-special} $\w(g_{2s})=\w(g_{2s-1})=\w(f_{2s-1})\neq m_{j_i}^{-1}$, hence $|(g_{2s-1}+g_{2s})(x_{2i}-x_{2i-1})|\leq 4 n_{2j-1}^{-2}$ by \ref{dep-seq-exact} and \ref{dep-seq-estimate}. Thus, using the estimate $\|x_i\|\leq 6$ for any $i\in F$ by \ref{dep-seq-exact}, we obtain 
\begin{align*}
\Big| f\Big(\sum_{i\in F}(x_{2i}-x_{2i-1}\Big)\Big| & \leq 5+|E(g_{2l-1}+g_{2l})(x_{2l}-x_{2l-1})|+|(g_{2r+1}+g_{2r+2})(x_{2r+2}-x_{2r+1})|\\
&+\sum_{i=l+1}^r|(g_{2i-1}+g_{2i})(x_{2i}-x_{2i-1})|
\\
&\leq 53+\sum_{i=l+1}^r|(g_{2i})(x_{2i})-(g_{2i-1})(x_{2i-1})|+\sum_{i=l+1}^r|(g_{2i})(x_{2i-1})|
\\
&+\sum_{i=l+1}^r|(g_{2i-1})(x_{2i})|
\end{align*}
Note that $(g_{2i})(x_{2i})=(g_{2i})(S^{k_i}x_{2i-1})=(R^{k_i}g_{2i})(x_{2i-1})=(g_{2i-1})(x_{2i-1})$,  $i=1,\dots,r$ by  \ref{dep-seq-special} and \ref{dep-seq-vectors-shifted}.

Assume $m\leqslant 0$. Then, by \ref{dep-seq-exact} and \ref{ss-2}, $\min\range (x_{2i})>\max\range(f_{2i-1})=\max\range(h_{2i-1})\geq \max\range( g_{2i-1})$, 
hence $(g_{2i-1})(x_{2i})=0$ for any $i=l+1,\dots,r$. By \ref{dep-seq-limit-intersection} there is at most one $s\in\{1,\dots, r\}$ with $2^m\range(g_{2s})\cap2^m\range(x_{2s-1})\subset \range(h_{2s})\cap2^m\range(x_{2s-1})=\range(f_{2s})\cap2^m\range(x_{2s-1})\neq \emptyset$, thus with $\range(g_{2s})\cap\range(x_{2s-1})\neq \emptyset$.

Assume $m>0$. As $\min\range(g_{2i})\geq\min\range(h_{2i})=\min\range (f_{2i})>\max\range(x_{2i-1})$ for any $i=1,\dots,r$ by \ref{dep-seq-exact} and \ref{ss-2}, $g_{2i}(x_{2i-1})=0$ for any $i=l+1,\dots,r$. By \ref{dep-seq-limit-intersection}, there is at most one $s\in\{1,\dots, r\}$ with $\range (g_{2s-1})\cap\range(x_{2s})\subset 2^m\neq \range(h_{2s-1})\cap \range(x_{2s})=2^m\range(f_{2s-1})\cap\range(x_{2s})\neq\emptyset$. 

Hence, for $m\in\Z$, we have
\[
\Big|f\Big(\sum_{i\in F}(x_{2i}-x_{2i-1})\Big)\Big|\leq 53 + 0+ 0+6=59,
\]
which ends the proof.
 \end{proof}

\subsection{Projections in $X$}  Applying results on dependent sequences from the previous section we prove an uncomplementability of $Y$ in $X$ (Theorem~\ref{not-complemented}) and uncomplementability of any subspace spanned by an unconditional basic sequence (Theorem~\ref{unc-not-complemented}) which combined with Propositions~\ref{isometric}, \ref{prop:odd unc} and Corollary \ref{reflexive} completes the proof of Theorem \ref{main}.

\begin{theorem}\label{not-complemented} The subspace $Y$ is not complemented in $X$.   
\end{theorem}

\begin{proof} Assume there is a bounded projection $\widetilde P:X\to Y$. Then, by the Krein-Milman-Rutman Theorem there is a block subspace $W$ with $W\cap Y=\{0\}$ and a bounded projection $P:W+Y\to Y$. For any $j\in \N$ choose, by Proposition~\ref{prop: choice of dep seq} with $N_i=1$ for each $i\in\N$, a $j$-dependent sequence 
$(x_i,f_i)_{i=1}^{n_{2j-1}}$ with
$(x_{2i-1})_{i=1}^{n_{2j-1}/2}\subset W$. By the definition of dependent sequences, 
$(x_{2i-1})_{i=1}^{n_{2j-1}/2}\subset Y$. Then for vectors
$$
y=\frac{m_{2j-1}}{n_{2j-1}}\sum_{i=1}^{n_{2j-1}/2}x_{2i}\in
Y, \ \ \ \
w=\frac{m_{2j-1}}{n_{2j-1}}\sum_{i=1}^{n_{2j-1}/2}x_{2i-1}\in W
$$
by Proposition~\ref{dep-est} we have $\|y-w\|\leq 240m_{2j-1}^{-1}$,
whereas
$$
\|y+w\|\geq \frac{1}{m_{2j-1}}\sum_{i=1}^{n_{2j-1}}f_i(w+z)\geq 1
$$
which, as $j$ is arbitrarily big, shows that $\inf\{\|y-w\|:\|y\|=1=\|w\|, y\in Y, w\in W\}=0$, contradicting the boundedness of $\widetilde P$. \end{proof}

\begin{theorem}\label{unc-not-complemented} The space $X$ admits no complemented subspaces with an  unconditional basis.
\end{theorem}
\begin{proof} Assume towards contradiction that some  normalized unconditional basic sequence $(\widetilde w_n)_{n=1}^\infty$ spans a complemented subspace of  $X$. By Rosenthal $\ell_1$ theorem and reflexivity of $X$ (Corollary~\ref{reflexive}), we can pass to a weakly null subsequence of $(w_n)_{n=1}^\infty$. By Krein-Milman-Rutman Theorem and unconditionality of   $(\widetilde w_n)_{n=1}^\infty$ there is a normalized  unconditional block sequence $(w_n)_{n=1}^\infty$, spanning a block subspace complemented in $X$ by a bounded projection $P$. Let $C\geqslant  0$ be the unconditional constant of $(w_n)_{n=1}^\infty$. 

Denote by $P_N:X\to \operatorname{span} \{e_1,\dots,e_N\}$, $N\in\N$, the canonical coordinate projections. For any $j\in \N$ we shall choose inductively integers $N_1<\dots<N_{n_{2j-1}}$, intervals $J_1<\dots<J_{n_{2j-1}/2}$ in $\N$ and a $j$-dependent sequence $(x_i,f_i)_{i=1}^{n_{2j-1}}$, satisfying the following for any $i=1,\dots,n_{2j-1}/2$:
\begin{enumerate}[(i)]
\item\label{uzu-not-1} $\supp w_n\subset (N_{2i-1},N_{2i})$, for each $n\in J_i$,   
\item\label{uzu-not-2}  $x_{2i-1}\in\spa \{w_n: n\in J_i\}$, $x_{2i}\in X_{N_{2i}}$, 
\item\label{uzu-not-3} $\|P_{r_i} P(Id-P_{N_{2i}})\|\leq (6C)^{-1}8^{-i}$, with $r_i=\max\supp w_{\max J_i}$,
\item\label{uzu-not-4} $\|(Id-P_{N_{2i+1}})Px_{2i}\|\leq C^{-1}8^{-i}$.
\end{enumerate}

Assuming objects above are chosen, let  $J=\bigcup_{i\in\N}J_i$. The subspace $\overline{\spa}\{w_n: n\in J\}$ is complemented in $X$ by the projection $Q_J P$, where $Q_J:\overline{\spa}\{w_n: n\in \N\}\to \overline{\spa}\{w_n: n\in J\}$ is the coordinate projection with $\|Q_J\|\leqslant C$. By \ref{uzu-not-1}-\ref{uzu-not-4}, and Definition~\ref{def-exact}\ref{def-exact:2}, 
\begin{equation}
\|Q_J Px_{2i}\|\leq \|Q_J P_{r_i} Px_{2i}\|+ \|Q_J (Id-P_{N_{2i+1}}) Px_{2i}\|\leq 4^{-i}, \quad i=1,\dots,n_{2j-1}/2    \label{eq}
\end{equation}
Then vectors
$$
w=\frac{m_{2j-1}}{n_{2j-1}}\sum_{i=1}^{n_{2j-1}/2}x_{2i-1}\in
\spa\{w_n: n\in J\}, \ \ \  x=\frac{m_{2j-1}}{n_{2j-1}}\sum_{i=1}^{n_{2j-1}/2}x_{2i}\in X, 
$$
satisfy $\|w-x\|\leq 240m_{2j-1}^{-1}$ by Proposition~\ref{dep-est}, whereas, by \eqref{eq},
$$
\|(Q_J P)(w-x)\|\geq\|w\|-\|(Q_J P)x\|\geq \frac{1}{m_{2j-1}}\sum_{i=1}^{n_{2j-1}}f_i(w) -\frac{m_{2j-1}}{2n_{2j-1}}\geq \frac{1}{4}
$$
which, as $j$ is arbitrarily big, contradicts the continuity of $Q_J P$. 

We choose objects satisfying \ref{uzu-not-1}-\ref{uzu-not-4} using Proposition~\ref{prop: choice of dep seq} and the following.
\begin{remark}\label{remark easy}
For any finite dimensional operator $T$ on $X$ and $\varepsilon>0$ there is $M\in\N$ with $\|T(Id - P_M)\|<\varepsilon$. Otherwise, by compactness argument, one can pick a bounded block sequence $(z_n)_{n=1}^\infty\subset X$ (hence weakly null) with $(Tz_n)_{n=1}^\infty$ converging in norm to a non-zero vector. 
\end{remark}
We shall choose inductively  $(N_i)_{i=1}^{n_{2j-1}}$, alongside with $(x_i,f_i)_{i=1}^{n_{2j-1}}$ according to Proposition~\ref{prop: choice of dep seq}, $(J_i)_{i=1}^{n_{2j-1}/2}$ chosen  by \ref{uzu-not-2} and $(r_i)_{i=1}^{n_{2j-1}/2}$ defined by \ref{uzu-not-3}. 

Let $N_1=1$, and take $x_1\in\spa\{w_n:n\in J_1\}$, provided by  Proposition~\ref{prop: choice of dep seq}. By Remark \ref{remark easy} for $T=P_{r_1}P$ and $\varepsilon=(48C)^{-1}$ pick $N_2>r_1$ satisfying \ref{uzu-not-3}.  Assuming that for $i<n_{2j-1}/2$, the first $2i$ elements are chosen, pick $N_{2i+1}>N_{2i}$ satisfying \ref{uzu-not-4}. Take $x_{2i+1}\in\spa\{w_n: n\in J_{i+1}\}$ provided by Proposition~\ref{prop: choice of dep seq}. By Remark~\ref{remark easy} for $T=P_{r_{i+1}}P$ and $\varepsilon=(6C)^{-1}8^{-i}$ take $N_{2i+2}>\max\{N_{2i+1},r_{i+1}\}$ satisfying \ref{uzu-not-3}, completing the inductive step. 
\end{proof}
\begin{proof}[Proof of Theorem~\ref{main}]    
Propositions~\ref{isometric} and \ref{prop:odd unc} combined with Corollary~\ref{reflexive}, Theorems~\ref{unc-not-complemented} and \ref{not-complemented} yield that the space $X$ satisfies the assertion of Theorem~\ref{main}. 
\end{proof}

\end{document}